\documentclass[11pt]{article}

\usepackage{amsmath}
\usepackage{amssymb}
\usepackage{graphicx}

\numberwithin{equation}{section}

\title{\bf Triangulated 3-Manifolds:
 from Haken's normal surfaces to Thurston's algebraic equation}
\author{Feng Luo\thanks{Partially Supported by a NSF grant}
}

\begin{document}
\maketitle \baselineskip.45cm



\begin{abstract}
\noindent  We give a brief summary of some of our work and our
joint work with  Stephan Tillmann on solving Thurston's equation
and Haken equation on triangulated 3-manifolds in this paper.
Several conjectures on the existence of solutions to Thurston's
equation and Haken equation are made. Resolutions of these
conjecture will lead to a new proof of the Poincar\'e conjecture
without using the Ricci flow. We approach these conjectures by a
finite dimensional variational principle so that its critical
points are related to solutions to Thurston's gluing equation and
Haken's normal surface equation. The action functional is the
volume. This is a generalization of an earlier program by Casson
and Rivin for compact 3-manifolds with torus boundary.
\end{abstract}

\section{Introduction}





\bigskip

This paper is based on several talks given by the author at the
conference  ``Interactions Between Hyperbolic Geometry, Quantum
Topology and Number Theory " at Columbia University in 2009 and a
few more places. The goal of the paper is to give a quick summary
of some of our work \cite{Lu1} and our joint work with Stephan
Tillmann \cite{LT1}, \cite{LT2} on triangulated 3-manifolds. Our
work is an attempt to connect geometry and topology of compact
3-manifolds from the point of view of triangulations. We will
recall Haken's normal surface theory, Thurston's work on
construction of hyperbolic structures, Neumann-Zagier's work, the
notion of angle structures introduced by Casson, Rivin and
Lackenby, and the work of several other people. One important
point we would like to emphasize is the role that the
Neumann-Zagier Poisson structure plays in these theories. It is
conceivable that the Neumann-Zagier Poisson structure will play an
important role in discretization and quantization of SL(2,$\bold
C$) Chern-Simons theory in dimension three.

A combination of the recent work of Segerman-Tillmann \cite{ST},
Futer-Gu\'eritaud \cite{FG}, Luo-Tillmann \cite{LT2} and
\cite{Lu1} has prompted us to make several conjectures on the
solutions of Thurston's equation and Haken's normal surface
equations. The resolution of some of these conjectures will
produce a new proof of the Poincar\'e conjecture without using the
Ricci flow method.

\medskip
Let us begin with a recall of closed triangulated pseudo
3-manifolds. Take a disjoint union of tetrahedra. Identify
codimension-1 faces of tetrahedra in pairs by affine
homeomorphisms. The quotient space is a \it triangulated closed
pseudo 3-manifold. \rm  (See \S 2.1 for more details). In
particular, closed triangulated 3-manifolds are closed
triangulated pseudo 3-manifolds  and ideally triangulated
3-manifolds are pseudo 3-manifolds with vertices removed. Given a
closed triangulated oriented pseudo 3-manifold, there are linear
and algebraic equations associated to the triangulation. Besides
the homology theories, the most prominent ones are Haken's
equation of normal surfaces \cite{Ha} and Thurston's algebraic
gluing equation for construction of hyperbolic metrics \cite{Thu1}
using hyperbolic ideal tetrahedra. Haken's theory is topological
and studies surfaces in 3-manifolds and Thurston's equation is
geometric and tries to construct hyperbolic metrics from the
triangulation. In the most general setting, Thurston's equation
tries to find representations of the fundamental group into
$PSL(2, \bold C)$ (\cite{Yo}). Much work has been done on both
normal surface theory and Thurston's equation with fantastic
consequences in the past fifty years.

Haken's normal surface equation is linear. A basis for the
solution space was found recently by Kang-Rubinstein \cite{KR}. In
particular, there are always solutions to Haken's equation with
non-zero quadrilateral coordinates.  The situation for solving
Thurston's equation is different. The main problem which motivates
our investigation is the following.

\medskip
\noindent {\bf Main Problem} \it Given a closed oriented
triangulated pseudo 3-manifold $(M, \bold T)$, when does there
exist a solution to Thurston's gluing equation? \rm

The most investigated cases in solving Thurston's equation  are
associated to ideal triangulated 3-manifolds with torus boundary
so that the complex numbers $z$ are in the upper-half plane  (see
for instance \cite{Thu1},  \cite{Til2}, \cite{Dun}, \cite{PW} and
many others). These solutions are closely related to the
hyperbolic structures. However, we intend to study Thurston's
equation and its solutions in the most general setting of closed
oriented triangulated pseudo 3-manifolds, in particular, on closed
triangulated 3-manifolds.
 Even though a solution to Thurston's equation in the
general setting does not necessarily produce a hyperbolic
structure, one can still obtain important information from it. For
instance, it
 was observed in \cite{Yo} (see also \cite{NZ},
\cite{ST}) that each solution of Thurston's equation produces a
representation of the fundamental group of the pseudo 3-manifold
with vertices of the triangulation removed to $PSL(2, \bold C)$. A
simplified version of a recent theorem of Segerman-Tillmann
\cite{ST} states that

\medskip
\noindent {\bf Theorem 1.1}(Segerman-Tillmann) \it If $(M, \bold
T)$ is a closed triangulated oriented 3-manifold so that the
triangulation supports a solution to Thurston's equation, then
each edge in $\bold T$ either has two distinct end points or  is
homotopically essential in $M$. \rm

\medskip
In particular, their theorem says any one-vertex triangulation of
a simply connected 3-manifold cannot support a solution to
Thurston's equation.  A combination of theorem 1.1 and a result of
\cite{Yang} gives an interesting solution to the main problem for
closed 3-manifold. Namely, a closed triangulated 3-manifold $(M,
\bold T)$ supports a solution to Thurston's equation if and only
if  there exists a representation $\rho: \pi_1(M) \to PSL(2, \bold
C)$ so that $\rho([e]) \neq 1$ for each edge $e$ having the same
end points. The drawback of this solution is that the
representation $\rho$ has to be a priori given.
\medskip

Our recent work \cite{Lu1} suggests another way to resolve the
main problem using Haken's normal surface equation.  To state the
corresponding conjecture, let us recall that a solution to Haken's
normal surface equation is said to be of \it 2-quad-type \rm if it
has exactly one or two non-zero quadrilateral coordinates. A \it
cluster of three 2-quad-type solutions \rm to Haken's equation
consists of three 2-quad-type solutions $x_1$, $x_2$ and $x_3$ so
that there is a tetrahedron containing three distinct
quadrilaterals $q_1, q_2, q_3$ with $x_i(q_i) \neq 0$ for
$i=1,2,3$. A triangulation of a 3-manifold is called \it minimal
\rm if it has the smallest number of tetrahedra among all
triangulations of the 3-manifold.

The main focus of our investigation will be around the following
conjecture. We thank Ben Burton and Henry Segerman for  providing
supporting data which helped us formulating it in the current
form.

\medskip
\noindent {\bf Conjecture 1} (Haken-Thurston Alternative)  \it For
any closed irreducible orientable minimally triangulated
3-manifold $(M, \bold T)$, one of the two holds:
\medskip

(1) there exists a solution to Thurston's equation associated to
the triangulation, or
\medskip

(2) there exists a cluster of three 2-quad-type solutions to
Haken's normal surface equation.   \rm

\medskip

Using a theorem of Futer-Gu\'eritaud, we proved in \cite{Lu1} the
following result which supports conjecture 1.

\medskip
\noindent {\bf Theorem 1.2} \it Suppose $(M, \bold T)$ is a closed
triangulated oriented pseudo 3-manifold. Then either there exists
a solution to the generalized Thurston equation or there exists a
cluster of three 2-quad-type solutions to Haken's normal surface
equation.  \rm

In our joint work with Tillmann \cite{LT2}, using
Jaco-Rubinstein's work \cite{LT2}, we proved the following theorem
concerning the topology of 3-manifolds satisfying part (2) of
conjecture 1.

\medskip
\noindent {\bf Theorem 1.3}(\cite{LT2}) \it Suppose $(M, T)$ is a
minimally triangulated orientable closed 3-manifold so that there
exists a cluster of three 2-quad-type solutions to Haken's
equation.  Then,

(a) $M$ is reducible, or

(b) $M$ is toroidal, or

(c) $M$ is a Seifert fibered space, or

(d) $M$ contains the connected sum $\#_{i=1}^3 RP^2$ of three
copies of the projective plane. \rm

\medskip
Using theorems  1.1 and 1.3,  one can deduce the Poincar\'e
conjecture from conjecture 1 (without using the Ricci flow) as
follows. Suppose $M$ is a simply connected closed 3-manifold. By
the Kneser-Milnor prime decomposition theorem, we may assume that
$M$ is irreducible. Take a minimal triangulation $\bold T$ of $M$.
By the work of Jaco-Rubinstein on 0-efficient triangulation
\cite{JR}, we may assume that $\bold T$ has only one vertex, i.e.,
each edge is a loop. By Segerman-Tillmann's theorem above, we see
that $(M, \bold T)$ cannot support a solution to Thurston's
equation. By conjecture 1, there exists a cluster of three
2-quad-type solutions to Haken's equation. By theorem 1.3, the
minimality of $\bold T$ and irreducibility of $M$, we conclude
that $M = \bold S^3$.

Theorem 1.2 is proved in  \cite{Lu1} where we proposed  a
variational principle associated to the triangulation to approach
conjecture 1. In this approach, 2-quad-type solutions to Haken's
equation arise naturally from non-smooth maximum points.  We
generalize the notion of angle structures introduced by Casson,
Lackenby \cite{Lac} and Rivin \cite{Ri1}
 (for ideally triangulated cusped 3-manifolds) to  the circle-valued
angle structure (or $\bold S^1$-angle structure or $SAS$ for
short) and its volume for any closed triangulated pseudo
3-manifold. It is essentially proved in \cite{LT1} and more
specifically in \cite{Lu1} that an $SAS$ exists on any closed
triangulated pseudo 3-manifold $(M, \bold T)$. The space
$SAS(\bold T)$ of all circle-valued angle structures on $(M, \bold
T)$ is shown to be a closed smooth manifold. Furthermore, each
circle-valued angle structure has a natural volume defined by the
Milnor-Lobachevsky function. The volume defines a continuous but
not necessarily smooth volume function $vol$ on the space
$SAS(\bold T)$. In particular, the volume function $vol$  achieves
a maximum point in $SAS(\bold T)$. The two conclusions in theorem
1.2 correspond to the maximum point being smooth or not for the
volume function.

More details of the results obtained so far and our approaches to
resolve the conjecture 1 will be discussed in  sections 4 and 5.
We remark that conjecture 1 itself is independent of the angle
structures and there are other ways to approach it.

\medskip

There are several interesting problems arising from the approach
taken here. For instance, how to relate the critical values of the
volume function on  $SAS(\bold T)$ with the Gromov norm of the
3-manifold. The Gromov norm of a closed 3-manifold is probably the
most important topological invariant for 3-manifolds. Yet its
computation is not easy. It seems highly likely that for a
triangulation without a cluster of three 2-quad-type solutions to
Haken's equation, the Gromov norm of the manifold (multiplied by
the volume of the regular ideal tetrahedron) is among the critical
values of the volume function on $SAS(\bold T)$. In our recent
work with Tillmann and  Yang \cite{LTY}, we have solved this
problem for closed hyperbolic manifolds. An affirmative resolution
of this problem for all 3-manifolds may provide insights which
help to resolve the Volume Conjecture for closed 3-manifolds.

\bigskip Futer and Gu\'eritaud have written a very nice paper \cite{Futer} on volume
and angle structures which is closely related to the material
covered in this paper.

We remark that this is not a survey paper on the subject of
triangulations of 3-manifolds. Important work in the field, in
particular the work of Jaco-Rubinstein \cite{JR} on efficient
triangulations of 3-manifolds, is not discussed in the paper.

\medskip
The paper is organized as follows.  In section 2, we will recall
the basic material on triangulations and Haken's normal surface
theory. In section 3, we discuss Neumann-Zagier's Poisson
structures and Thurston's gluing equation. In section 4, we
discuss circle valued angle structures, their volume, some of our
work and a theorem of Futer-Gu\'eritaud.  In section 5, we
introduce a $\bold Z_2$ version of Thurston's equation ($\bold
Z_2$-taut structure).

\medskip
\noindent {\bf Acknowledgement.}  The work is supported in part by
the NSF. We thank the editors of the conference proceedings for
inviting  us to write the paper and S. Tillmann and the referee
for suggestions on improving the writing of this paper. We would
like to thank in particular David Futer and Fran{c}ois Gu\'eritaud
for allowing us to present their unpublished theorem. The proof of
this theorem was also supplied by them.









\section{\bf Triangulations and normal surfaces}

 The  normal surface theory, developed by Haken in the 1950's, is a
beautiful chapter in 3-manifold topology. In the late 1970's,
Thurston introduced the notion of spun normal surfaces and used it
to study 3-manifolds.

We will revisit the normal surface theory and follow the
expositions in \cite{JT} and \cite{Til} closely in this section.
Some of the notations used in this section are new. The work of
Tollefson, Kang-Rubinstein, Tillmann, and Jaco on characterizing
the quadrilateral coordinates of normal surfaces will be
discussed.

\subsection{\bf Some useful facts about tetrahedra}

The following lemma will be used frequently in the sequel. The
proof is very simple and will be omitted.  To start, suppose
$\sigma =[v_1, .., v_4]$ is a tetrahedron with vertices $v_1, ...,
v_4$ and edges $e_{ij}=\{v_i, v_j\} $, $i \neq j$. We call
$e_{kl}$ the \it opposite edge \rm of $e_{ij}$ if
$\{i,j,k,l\}=\{1,2,3,4\}$.

\medskip

\noindent {\bf Lemma 2.1.} \it Given a tetrahedron $\sigma$,
assign to each edge $e_{ij}$ a real number $a_{ij} \in \bold R$,
called the weight of $e_{ij}$.  Assume $\{i,j,k,l\}=\{1,2,3,4\}$.

(a)  If the sum of weights of opposite edges is a constant, i.e.,
$a_{ij}+ a_{kl}$ is independent of indices, then there exist real
numbers $b_1, .., b_4$ (weights at vertices) so that $$a_{ij} =
b_i + b_j.$$

(b) If  the sum of  weights of the edges from each vertex is a
constant, i.e., $a_{ij}+ a_{ik} + a_{il}$ is independent of
 indices, then weights of opposite edges are the
same, i.e.,
$$ a_{ij} = a_{kl}.$$

(c) If the tetrahedron $\sigma$ is oriented and edges are labelled
by $a,b,c$ so that opposite edges are labelled by the same letter
(see figure 1(a)), then the cyclic order $a \to b \to c \to a$ is
independent of the choice of the vertices and depends only on the
orientation of $\sigma$. \rm

\subsection{ Triangulated closed pseudo 3-manifolds and Haken's normal surface equation}

\medskip
Let $X$ be a union of finitely many disjoint oriented Euclidean
tetrahedra. The collection of all faces of tetrahedra in $X$ is a
simplicial complex $\bold T^*$ which is a triangulation of $X$.
Identify codimension-1 faces in $X$ in pairs  by affine
orientation-reversing homeomorphisms. The quotient space $M$ is a
closed oriented pseudo 3-manifold with a triangulation $\bold T$
whose simplices are the quotients of simplices in $\bold T^*$. Let
$V, E, F, T$ (and $V^*, E^*, F^*$ and $T^*$) be the sets of all
vertices, edges, triangles and tetrahedra in $\bold T$ (in $\bold
T^*$ respectively). The quotient of a simplex $x \in \bold T^*$
will be denoted by $[x]$ in $\bold T$. We call $x \in \bold T^*$
the unidentified simplex and $[x]$ the quotient simplex. Since the
sets of tetrahedra in $T^*$ and $T$ are bijective under the
quotient map, we will identify a tetrahedron $\sigma \in T^*$ with
its quotient $[\sigma]$, i.e., $\sigma =[\sigma]$ and $T=T^*$.

If $x, y \in V \cup E \cup F \cup T$  (or in $\bold T^*$), we use
$ x
> y$ to denote that $y$ is a face of $x$. We use $|Y|$ to denote
the cardinality of a set $Y$.

Note that in this definition of triangulation, we do not assume
that simplices in $\bold T$ are embedded in $M$. For instance, it
may well be that $|V| =1$. Furthermore, the non-manifold points in
$M$ are contained in  the set of vertices.

\medskip
According to Haken, a normal surface in a triangulated pseudo
3-manifold $M$ is an embedded surface $S \subset M$ so that for
each tetrahedron $\sigma$, topologically the intersection $S \cap
\sigma$ consists of a collection of planar quadrilaterals and
planar triangles, i.e., inside each tetrahedron, topologically the
surface $S$ looks like planes cutting through the tetrahedron
generically. Haken's theory puts this geometric observation into
an algebraic setting. According to \cite{Ha}, a \it normal arc \rm
in $X$ is an embedded arc in a triangle face so that its end
points are in different edges and a \it normal disk \rm in $X$ is
an embedded disk in a tetrahedron so that its boundary consists of
3 or 4 normal arcs. These are called \it normal triangles \rm and
\it normal quadrilaterals \rm respectively.    A normal isotopy is
an isotopy of $X$ leaving each simplex invariant. Haken's normal
surface theory deals with normal isotopy classes of normal disks
and normal surfaces. For simplicity, we will interchange the use
of normal disk with the normal isotopy class of a normal disk.

\begin{figure}[ht!]
\centering
\includegraphics[scale=0.7]{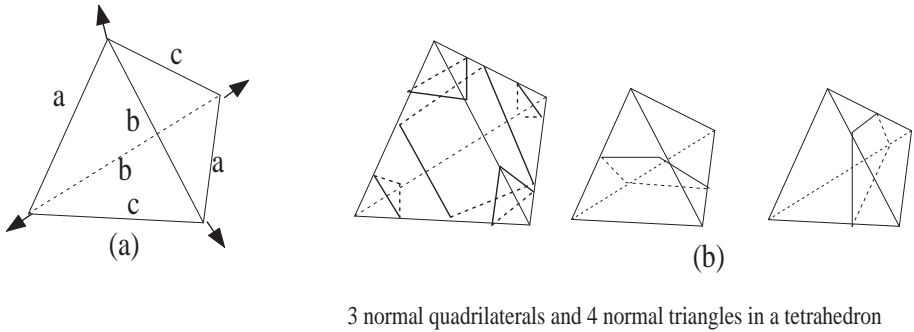}
\caption{$t, t'$ are normal triangles and $q, q'$ are normal
quadrilaterals} \label{figure 2.1}
\end{figure}






The projections of normal arcs and normal disks from $X$ to $M$
constitute normal arcs and normal disks in the triangulated space
$(M, \bold T)$.  For each tetrahedron, there are four normal
triangles and three normal quadrilaterals inside it up to normal
isotopy. See figure 1(b).  Note that there is a natural one-one
correspondence between normal disks in $T^*$ and $T$. In the
sequel, we will not distinguish normal disks in $\bold T$ or
$\bold T^*$ and we will use $\triangle$, $\Box$ to denote the sets
of all normal isotopy classes of normal triangles and
quadrilaterals in the triangulation $\bold T$ and also $\bold
T^*$. The set of normal arcs in $\bold T^*$ and $\bold T$ are
denoted by $\bold A^*$ and $\bold A$ respectively.

There are relationships among the sets $V, E, F, T, \triangle,
\Box, \bold A$. These incidence relations, which will be recalled
below, are the basic ingredients for defining Haken's and
Thurston's equations.

Take $ t \in \triangle$, $a \in \bold A$, $q \in \Box$, and
$\sigma \in T$.  The following notations will be used. We use $ a<
t$ (and $ a < q$) if there exist representatives $x \in a$, $y \in
t$ (and $z \in q$) so that $x$ is an edge of $y$ (and $z$). We use
$t \subset \sigma$ and $q \subset \sigma$ to denote that
representatives of $t$ and $q$ are in the tetrahedron $\sigma$. In
this case, we say the tetrahedron $\sigma$ \it contains \rm $t$
and $q$.

As a convention,  we will always use the letters $\sigma$, $e$ and
$q$ to denote a tetrahedron, an edge and a quadrilateral in the
triangulation $\bold T$ respectively.

The normal surface equation is a system of linear equations
defined in the space $\bold R^{\triangle} \times \bold R^{\Box}$,
introduced by W. Haken \cite{Ha}.  It is defined as follows. For
each normal arc $ a \in \bold A$, suppose $\sigma, \sigma'$ are
the two tetrahedra adjacent to the triangular face which contains
$a$. (Note that $\sigma$ may be $\sigma'$.)  Then there is a
homogeneous linear equation for $ x \in \bold R^{\triangle}\times
\bold R^{\Box}$ associated to $a$:

\begin{equation} \label{2.1}
 x (t) + x(q) = x(q') + x(t')
 \end{equation}
 where  $t, q \subset \sigma$, $t', q' \subset \sigma'$ and $t, t',
 q, q' > a$. See figure 2(a).


Recall that we identify the set of edges $E$ with the quotient of
$E^*$, i.e., $E =\{ [y] | y \in E^*\}$ where $[y] =\{ y' \in E^* |
y \sim y'\}$. The index $i: E^* \times \Box \to \bold Z$ is
defined as follows: $i(y, q) = 1$ if $y, q$ lie in the same
tetrahedron $\sigma \in \bold T^*$ so that $y \cap q =\emptyset$,
and $i(y, q) =0$ in all other cases. The index $i: E \times \Box
\to \bold Z$ is defined to be $i(e, q) =\sum_{ y \in e} i(y, q)$.
See figure 2(b) for a picture of $i(e,q)=1, 2$. For simplicial
triangulations, $i(e, q)=1$ means that the quadrilateral $q$ faces
the edge $e$ in a tetrahedron, i.e., $ q \cap e = \emptyset$ and
$e, q \subset \sigma$. In general, $i(e,q) \in \{0,1,2\}$.
However, for simplicial triangulations, $i(e,q)=0, 1$.

\begin{figure}[ht!]
\centering
\includegraphics[scale=0.7]{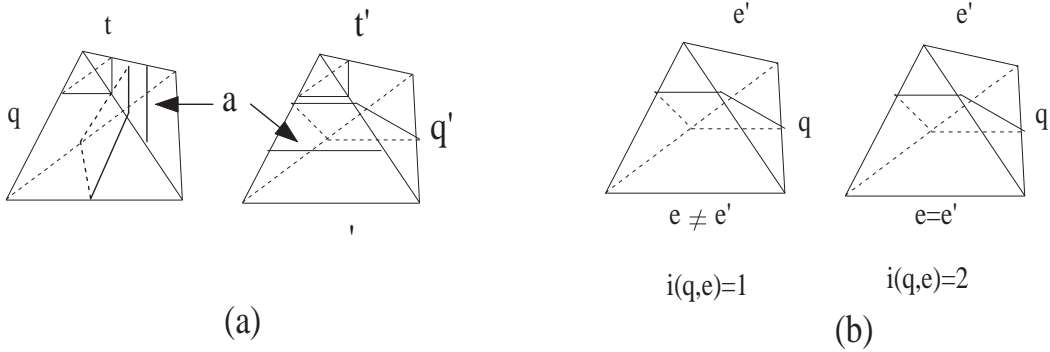}
\caption{incident indices     } \label{1.1: figure 1.1}
\end{figure}

\subsection{Normal surfaces and tangential angle structures}

Given  $x \in \bold R^{\triangle} \times \bold R^{\Box}$,  we will
call $x(t)$ ($t \in \Delta$) and $x(q)$ the $t$-coordinate and
$q$-coordinate (triangle and quadrilateral coordinates) of $x$.
Haken's normal surface equation addresses the following question.
Given a finite set of normal triangles and normal quadrilaterals
in a triangulation $\bold T$, when can one construct a normal
surface with these given triangles and quadrilaterals as its
intersections with the tetrahedra? Haken's equation (2.1) is a set
of necessary conditions. Spun normal surface theory addresses the
following question, first investigated by Thurston \cite{Thu1}.
Suppose we are given a finite set of quadrilaterals in each
tetrahedron. When can one construct a normal surface whose
quadrilateral set is the given one?  We can phrase it in terms of
the normal coordinates as follows. Given a vector $z \in \bold
R^{\Box}$, when does there exist a solution to Haken's equation
(2.1) whose projection to $\bold R^{\square}$ is $z$? The question
was completely solved in \cite{Tol}, \cite{KR}, \cite{Til} and
\cite{Jac}.  We will interpret their results in terms of angle
structures.

\bigskip
\noindent {\bf Definition 2.1.} A \it tangential angle structure
\rm on a triangulated pseudo 3-manifold $(M, \bold T)$ is a vector
$ x \in \bold R^{\Box}$ so that,

 for each tetrahedron $\sigma \in T$, \begin{equation} \sum_{ q
\in \Box, q \subset \sigma} x(q)=0, \end{equation}

and for each edge $e \in E$, \begin{equation} \sum_{q \in \Box}
i(e,q) x(q) =0. \end{equation}

The linear space of all tangential angle structures on $(M, \bold
T)$ is denoted by $TAS(\bold T)$ or $TAS$.
\medskip

Recall that a (Euclidean type) angle structure, introduced by
Casson, Rivin \cite{Ri1} and Lackenby \cite{Lac}, is a vector $ x
\in \bold R_{>0}^{\Box}$ so that for each tetrahedron $\sigma \in
T$,
\begin{equation} \sum_{ q \in \Box, q \subset \sigma}
x(q)=\pi, \end{equation} and
 for each $e \in E$, \begin{equation} \sum_{q \in \Box} i(e,q)
x(q) =2\pi. \end{equation}

These two conditions (2.4) and (2.5) have very natural geometric
meaning. Suppose a hyperbolic manifold admits a geometric
triangulation by ideal hyperbolic tetrahedra. The first equation
(2.4) says that a normal triangle in a hyperbolic ideal
tetrahedron is Euclidean and the second equation (2.5) says that
the sum of the dihedral angles around each edge is $2\pi$. By
definition, a tangential angle structure is a tangent vector to
the space of all angle structures.

The following is a result proved by Tollefson (for closed
3-manifolds), Kang-Rubinstein and Tillmann for all cases. The
result was also known to Jaco \cite{Jac}.  Let $\bf S_{ns}$ be the
space of all solutions to Haken's homogeneous linear equations
(2.1).

Given a finite set $X$, the standard basis of $\bold R^X$ will be
denoted by $X^*$ =\{$ x^* \in \bold R^X|  x \in X\}$ so that
$x^*(t) =0$ if $t \in  X -\{x\}$ and $x^*(x)=1$. We give $\bold
R^X$ the standard inner product $(\ ,\ )$ so that $X^*$ forms an
orthonormal basis.

\medskip
\noindent {\bf Theorem 2.2.} (\cite{Tol}, \cite{KR}, \cite{Til})
\it For a triangulated closed pseudo 3-manifold $(M, \bold T)$,
let $Proj_{\Box}: \bold R^{\triangle} \times \bold R^{\Box} \to
\bold R^{\Box}$ be the projection. Then
\begin{equation}
Proj_{\Box}(\bold S_{ns}) =TAS(\bold T)^{\perp},
\end{equation}
where $\bold R^{\Box}$ has the standard inner product so that
\{$q^* | q \in \Box\}$ is an orthonormal basis. \rm

\medskip
For a short proof of this theorem, see \cite{Lu1}. This result is
very important for us to relate normal surfaces to critical points
of the volume function on the space of all circle-valued angle
structures.

\section{\bf Neumann-Zagier Poisson structure and Thurston's gluing equation}

The Neumann-Zagier Poisson structure on $\bold R^{\Box}$,
introduced in \cite{NZ}, is of fundamental importance for studying
triangulated 3-manifolds and in particular for Thurston's gluing
equation.  We will recall its definition and derive some of its
properties in this section. See also \cite{Cho1} and \cite{Cho2}
for different proofs.

\subsection{The Neumann-Zagier Poisson structure}

Recall that our triangulated pseudo 3-manifolds $(M, \bold T)$ are
oriented so that each tetrahedron has the induced orientation.
Since a pair of opposite edges $\{e, e'\}$ in a tetrahedron
$\sigma$ is the same as a normal quadrilateral $q \subset \sigma$
with $i(e, q) \neq 0$, by lemma 2.1, for each tetrahedron $\sigma$
in $\bold T$, there exists a natural cyclic order on the three
quadrilaterals $q_1, q_2, q_3$ in $\sigma$. We denote the cyclic
order by $q_1 \to q_2 \to q_3 \to q_1$, and write $q \to q'$ in
$\Box$ if $q, q'$ are in the same tetrahedron and $ q \to q'$ in
the cyclic order. Define a map $w: \Box \times \Box \to \bold R$
by $w(q, q') =1$ if $q \to q'$, $w(q, q')=-1$ if $q' \to q$ and
$w(x,y)=0$ otherwise. The Neumann-Zagier skew symmetric bilinear
form, still denoted by $w: \bold R^{\Box} \times \bold R^{\Box}
\to \bold R$, is defined to be:

$$  w(x, y) = \sum_{q, q' \in \Box}  w(q,q') x(q) y(q').$$
From the definition, it is evident that $w(x, y)=-w(y,x)$.

The following was proved in \cite{NZ},

\noindent

\medskip
\noindent {\bf Proposition 3.1 (Neumann-Zagier)}. \it Suppose $(M,
\bold T)$ is a triangulated, oriented closed pseudo 3-manifold.
Then

(a) for any $q' \in \Box$, $\sum_{q \in \Box} w(q, q') =0$,

(b) for any pair of edges $e, e' \in E$,
$$ \sum_{ q, q' \in \Box}  i(e,q) i(e',q') w(q, q') =0.$$
\rm

Let $Z$ be the linear subspace $\{ x \in \bold R^{\Box} |$ for all
$\sigma \in T$, $\sum_{ q \subset \sigma} x(q) =0$\}. Then the \it
Neumann-Zagier symplectic 2-form \rm is the restriction of $w$ to
$Z^2$. It provides an identification between $Z$ and the dual
space $Z^*$.

A simple property of the Neumann-Zagier form is the following
identity. For any $q_1, q_2 \in \Box$,
\begin{equation}
\sum_{ q \in \Box} w(q_1, q) w(q, q_2) = \left\{
     \begin{array}{lr}

0, &    q_1, q_2  \quad \text{not in a tetrahedron}\\
-2, &  q_1=q_2  \\
1, &  q_1 \neq q_2 \quad \text{and}  \quad q_1, q_2 \subset \sigma

     \end{array}
   \right.
\end{equation}
If $y \in Z$, then
\begin{equation}
\sum_{q, q_2 \in \Box} w(q_1, q) w(q, q_2) y(q_2) = -3 y(q_1).
\end{equation}

Indeed,  by (3.1), the left-hand-side of (3.2) is equal to
$-2y(q_1)+ y(q_3)+y(q_4)$ where $q_1, q_3, q_4$ are three
quadrilaterals in a tetrahedron. Since $y(q_1)+ y(q_3)+ y(q_4)=0$
by definition of $Z$, equation (3.2) follows.  For any $q' \in
\Box$, the vector $y(q)= w(q, q')$ is an element in $Z$ by
proposition 3.1(a). Putting this $y=y(q_2)= w(q_2, q_4)$ into
identity (3.2), we obtain, for any $q_1, q_4 \in \Box$,

\begin{equation}
\sum_{q_2, q_3 \in \Box} w(q_1, q_2) w(q_2, q_3) w(q_3, q_4) = -3
w(q_1, q_4).
\end{equation}

 We will identify the dual space $(\bold R^X)^*$ with $\bold
R^X$ via the standard inner product $( \  ,  \ )$ where $X^*$ is
an orthonormal basis.

For a triangulated pseudo 3-manifold $(M, \bold T)$, define the
linear  map $A: Z \to \bold R^E$ by

\begin{equation}
A(x)(e) = \sum_{ q} i(e, q) x(q). \end{equation}  Note that the
space of all tangential angle structures $TAS$ is exactly equal to
$ker(A)$.

\medskip
\noindent {\bf Lemma 3.2. } \it Suppose $(M, \bold T)$ is
oriented. The dual map $A^*: \bold R^E \to Z$, where the dual
spaces of $\bold R^E$ and $Z$ are identified with themselves via
the standard inner product $(,)$ on $\bold R^E$ and $w$ on $Z$, is
$$A^*(x)(q) =\frac{1}{3}\sum_{e \in E} W(e, q) x(e),$$ where
$$W(e, q) = \sum_{q' \in \Box} i(e,q') w(q', q).$$ \rm

\medskip \noindent {\bf Proof.}  We need to show for any $x \in
\bold R^E$ and $y \in Z$, $$ (A(y), x) = w(y, A^*(x)).$$

Indeed, the left-hand-side of it is
$$ \sum_{e} A(y)(e) x(e) = \sum_{e, q} x(e) i(e, q) y(q).$$

The right-hand-side of it is
$$\sum_{q', q'' \in \Box} y(q') w(q', q'') A^*(x) (q'')$$
$$=\frac{1}{3} \sum_{q', q'', e} y(q') w(q', q'') W(e, q'') x(e)$$
$$ = \frac{1}{3} \sum_{q', q'', e, q} y(q') w(q', q'') i(e, q) w(q, q'') x(e)$$

$$=-\frac{1}{3} \sum_{e, q} i(e, q) x(e)  \sum_{ q', q''} y(q') w(q', q'') w(q'', q)$$
$$ =\sum_{e, q} i(e, q) x(e) y(q).$$

Here the last equation comes from (3.3).  This ends the proof.

\medskip
 Let $B: \bold R^E \to \bold R^V$ be the map $B(x)(v) =
\sum_{ e > v} x(e)$. If both end points of $e$ are $v$, then the
edge $e$ is counted twice in the summation $\sum_{ e > v} x(e)$.
 The dual map $B^*: \bold R^V
\to \bold R^E$ is given by $ B^*(y)(e) = \sum_{ v < e} y(v).$

\medskip
\noindent {\bf Theorem 3.3 (Neumann-Zagier)} \it For any oriented
triangulated closed pseudo 3-manifold $(M, \bold T)$, the
sequences

\begin{equation*}
Z \xrightarrow{A} \bold R^E \xrightarrow{B} \bold R^V
\xrightarrow{} 0
\end{equation*}
and
\begin{equation*}
0 \xrightarrow{} \bold R^V \xrightarrow{B^*} \bold R^E
\xrightarrow{A^*} Z
\end{equation*}
are exact.  Furthermore, if $x, y \in \bold R^E$, then $$ w(
A^*(x), A^*(y)) =0.$$

\rm

\medskip
\noindent {\bf Proof.} (See also \cite{Cho1}.) Since the second
sequence is the dual of the first, it suffices to prove that one
of them is exact. First, $BA=0$ follows from the definition of
$Z$. Furthermore, it is easy to see that $B^*$ is injective.
Indeed, if $B^*(y) =0$ for some $y \in \bold R^V$, then by
definition, $y(v)+ y(v')=0$ whenever $v, v'$ form the end points
of an edge. Now for any $v \in V$, take a triangle in $\bold T$
with vertices $v_1=v, v_2,$ and $v_3$. Then
 equations $y(v_i) + y(v_j) =0$ for $ i \neq j$ in $\{1,2,3\}$
 imply that $y(v_i)=0$, i.e., $y(v)=0$. It remains to prove that $ker(A^*) \subset
 Im(B^*)$. Suppose $x \in \bold R^E$ so that $A^*(x) =0$, i.e.,
 for all $q \in \Box$, $$ A^*(x)(q)  = \frac{1}{3} \sum_{ e} W(e, q) x(e)=0.$$
 Spelling out the details of the above equation, we see that it is
 equivalent to

 $$ x(e_1) + x(e_2) = x(e_3) + x(e_4)$$ whenever
 $\{e_1, e_2\}$ and $\{e_3, e_4\}$ are two pairs of opposite edges
 in a tetrahedron $\sigma$ in $\bold T$. Fix a tetrahedron $\sigma$, and
 consider $x(e)$ as weights on the edges of $\sigma$. By lemma
 2.1, there exists a map $y: \{ (v, \sigma) | v < \sigma, v \in
 V\} \to \bold R$ so that

\begin{equation}
  x(e) = \sum_{v < e}  y(v, \sigma).
\end{equation}

 We claim that the above equation implies that $y(v, \sigma) =
 y(v, \sigma')$ for any other tetrahedron $\sigma' > v$. Assuming
 this claim, and taking $y(v)= y(v, \sigma)$, then we have $x(e) =\sum_{v < e} y(v)$, i.e., $x
 =B^*(y)$, or $ x \in Im(B^*)$.

 To see the claim, let us first assume that $\sigma$ and $\sigma'$
 share a common triangle face which has $v$ as a vertex. Say the
 three vertices of the triangle face are $v_1=v, v_2,$ and $v_3$.
 Then equation (3.4) says that
 $$ y(v_i, \sigma) + y(v_j, \sigma) = y(v_i, \sigma') + y(v_j,
 \sigma'),$$ for $ i \neq j $ in $\{1,2,3\}$. The common value is
 $x_{ij}=x(\{v_i, v_j\})$. This system of three equations has a unique
 solution, namely $y(v_i, \sigma) = y(v_i, \sigma') = ( x_{ik}+
 x_{ij} - x_{jk})/2$ for $\{i,j,k\}=\{1,2,3\}$. Now in general, if
 $\sigma$ and $\sigma'$ are two tetrahedra in $\bold T$ which have
 a common vertex $v$,  by the definition of pseudo 3-manifolds,
 there exists a sequence of tetrahedra $\sigma_1 =\sigma,
 \sigma_2, ..., \sigma_n=\sigma'$ so that for each index $i$,
 $\sigma_i, \sigma_{i+1}$ share a common triangle face which has
 $v$ as a vertex.  Thus, by repeating the same argument just given, we
 see that $y(v, \sigma) = y(v, \sigma').$

To see the last identity,

$$ w(A^*(x), A^*(y))
=\sum_{q_1, q_2} w(q_1, q_2) A^*(x)(q_1) A^*(y)(q_2)$$
$$=\frac{1}{9}\sum_{q_1, q_2, e_1, e_2} w(q_1, q_2) W(e_1, q_1) x(e_1) W(e_2,
q_2) y(e_2)$$
$$=\frac{1}{9} \sum_{q_1, q_2, q_3, q_4, e_1, e_2} w(q_1, q_2) i(e_1, q_3)
w(q_3, q_1) i(e_2, q_4) w(q_4, q_2) x(e_1) y(e_2)$$

$$=-\frac{1}{9}\sum_{q_3, q_4, e_1, e_2} i(e_1, q_3) i(e_2, q_4) x(e_1)
y(e_2) \sum_{q_1, q_2} w(q_3, q_1) w(q_1, q_2) w(q_2, q_4). $$

By (3.3) and proposition 3.1(b),  the above is equal to
$$ = \frac{1}{3} \sum_{q_3, q_4, e_1, e_2} i(e_1, q_3) i(e_2, q_4) x(e_1)
y(e_2) w(q_3, q_4) =0.$$ This ends the proof.

\subsection{\bf Thurston's equation}

Let us recall briefly  Thurston's gluing equation \cite{Thu1} on a
triangulated closed oriented pseudo 3-manifold $(M, \bold T)$.
Assign each edge in each tetrahedron in the triangulation $\bold
T$ a complex number $z \in \bold C-\{0,1\}$. The assignment is
said to satisfy the \it generalized Thurston algebraic equation
\rm if

\medskip

\noindent (a) opposite edges of each tetrahedron have the same
assignment;

\noindent (b) the three complex numbers assigned to three pairs of
opposite edges in each tetrahedron are  $z$, $\frac{1}{1-z}$ and
$\frac{z-1}{z}$ subject to an orientation convention; and

\noindent (c) for each edge $e$ in the triangulation, if \{$z_1,
..., z_k$\} is the set of all complex numbers assigned to the edge
$e$ in the various tetrahedra adjacent to $e$,
 then
\begin{equation}  \prod_{i=1}^k z_i = \pm 1.
\end{equation}

If the right-hand-side of (3.6) equals 1 for all edges, we say
that the assignment satisfies \it Thurston algebraic equation. \rm

Since a  pair of opposite edges  in a tetrahedron  is the same as
the normal isotopy class of a quadrilateral, we see that
Thurston's equation is defined on $\bold C^{\Box}$. To be more
precise, given
 $z \in \bold C^{\Box}$, we say $z$ satisfies the generalized
Thurston equation, if the following assertions are satisfied:

(1) if $q \to q'$ in $\Box$, then $z(q') = \frac{1}{1-z(q)}$, and

(2) if $e \in E$, then

\begin{equation}
\prod_{ q}  z(q) ^{i(e, q)} = \pm 1.
\end{equation}

If the right-hand-side of (3.7) equals $1$ for all edges, we say
$z$ satisfies Thurston's equation.

This equation was introduced by Thurston in \cite{Thu1} in 1978.
He used it to construct the complete hyperbolic metric on the
figure-eight knot complement. Since then, many authors have
studied  Thurston's equation. See for instance \cite{NZ},
\cite{petro}, \cite{Cho1}, \cite{hodg}, \cite{Yo}, \cite{Til2} and
others. This equation was originally defined for ideal
triangulated 3-manifolds with torus boundary, i.e., closed
triangulated pseudo 3-manifolds $(M, \bold T)$ so that each vertex
link is a torus. We would like to point out that Thurston's
equation (3.6) is defined on any closed triangulated oriented
pseudo 3-manifold. It was first observed by Yoshida \cite{Yo},  a
solution to Thurston's equation produces a representation of the
fundamental group $\pi_1(M -\bold T^{(0)})$ to $PSL(2, \bold C)$
where $\bold T^{(0)}$ is the set of all vertices. Thus, in the
broader setting, solving Thurston's equation amounts to find
$PSL(2, \bold C)$ representations of the fundamental group. The
recent work of \cite{kk} seems to have rediscovered equation (3.6)
independently while working on TQFT.

Let $\bold D(\bold T)$ be the space of all solutions to Thurston's
equation defined in $\bold C^{\Box}$. By definition, $\bold
D(\bold T)$ is an algebraic set. There are several very nice
results known for $\bold D(\bold T)$.  Let $\bold H =\{ w \in
\bold C| \quad im(w)
>0\}$ be the upper-half-plane.

\medskip
\noindent {\bf Theorem 3.4 (Choi \cite{Cho1})}. \it The set $\bold
D(\bold T) \cap \bold H^{\Box}$ is a smooth complex manifold. \rm

\medskip
Her proof makes an essential use of Neumann-Zagier's symplectic
form (theorem 3.3).

Another result on Thurston's equation is in the work of Tillmann
\cite{Til2} and Yoshida \cite{Yo} relating degenerations of
solutions of Thurston's equation to normal surface theory. See
also the work of \cite{Ka} and \cite{Se}. The geometry behind
their construction was first observed by Thurston \cite{Thu2}.
Though this work does not address conjecture 1 in the
introduction, it does indicate  a relationship between Thurston's
equation and Haken's equation.

Here is Tillmann's construction. Suppose $z_n \in \bold D(\bold
T)$ is an unbounded sequence of solutions to Thurston's equation
(3.7) so that for each $q \in \Box$,

$$  u(q) = \lim_{ n \to \infty}  \frac{\ln| z_n (q)|}{\sqrt{1+ \sum_{q' \in \Box}
(\ln|z_n(q')|)^2}}$$ exists in $\bold R$.

Take the logarithm of equation (3.7) for $z_n$, divide the
resulting equation by $\sqrt{1+\sum_{q' \in \Box}( \ln|
z_n(q')|})^2$, and let $n \to \infty$. We obtain, for each edge $e
\in E$,
\begin{equation}
\sum_{q} i(e, q) u(q) =0. \end{equation}

By definition, $u(q)=0$ unless $\lim_{n \to \infty} z_n(q) =0, $
or $\infty$. Furthermore, if $\lim_{n} z_n(q) =1$ and $q' \to q
\to q''$, then $\lim_{ n} z_n(q'') =  \lim_n \frac{1}{1-z_n(q)}
=\infty$ and $\lim_{n} z_n(q') = \lim_n \frac{z_n(q)-1}{z_n(q)}=
0$ so that $\lim_n z_n(q') z_n(q'')=-1$. This implies that
$u(q')+u(q'')=0$ and $u(q'') \geq 0$. Let
 $I =\{ q \in \Box | \lim_{ n} z_n(q) =1\}$ and for $q \in I$, let  $a_q = u(q'') \geq 0$ where $q \to q''$. Then
$$ u =\sum_{ q \in I} a_q  \sum_{q'}  w(q, q') (q')^* \in \bold R.$$
Substitute into (3.8), we obtain for each $e \in E$,
\begin{equation}
 \sum_{ q \in I} v(q) W(e, q) =0
\end{equation}
where $v=\sum_{q \in I} a_q q^*$. Equation (3.9) appeared in the
work of Tollefson \cite{Tol} in which he proved that, if $(M,
\bold T)$ is a closed 3-manifold, then (3.9) gives a complete
characterization of the quadrilateral coordinates of solutions to
Haken's equation. Namely, if $M$ is closed,  a vector $v \in \bold
R^{\Box}$ is in $Proj_{\Box}(\bold S_{ns})$ if and only if (3.9)
holds for all $e \in E$.   Thus, by Tollefson's theorem, the
specific $ v= \sum_{q \in I} a_q q^* $ belongs to
$Proj_{\Box}(\bold S_{ne})$. As a consequence, one has,

\medskip
\noindent {\bf Theorem 3.5 }(Tillmann) \it For a closed
triangulated 3-manifold $(M, \bold T)$, the logarithmic limits of
$\bold D(\bold T)$ correspond to solutions of Haken's normal
surface equation. \rm

\medskip
We remark that Tillmann's theorem in \cite{Til2} is more general
and works for all pseudo 3-manifolds. We state it in the above
form for simplicity.    Furthermore, Tillmann observed in
\cite{Til2} that the solution $v$ has the property that there is
at most one non-zero quadrilateral coordinate in each tetrahedron.
Thus if all coefficients $a_q$ are non-negative integers, then the
vector $v$ produces an embedded normal surface in the manifold.

It follows from the definition that for each $e \in E$, the vector
\begin{equation}
 u_e = \sum_{q} w(e, q)q^*
\end{equation}
 is in $TAS(\bold T)$. What Tollefson proved, using the
language of TAS, is that for a closed triangulated 3-manifold $(M,
\bold T)$, the set $\{ u_e | e \in E\}$ generates the linear space
$TAS(\bold T)$. A generating set for $TAS(\bold T)$ for all closed
triangulated pseudo 3-manifolds $(M, \bold T)$ was found in the
work of Kang-Rubinstein \cite{KR} and Tillmann \cite{Til}.

In the recent work \cite{Yang}, Yang is able to construct many
solutions of  Thurston's equation on closed triangulated
3-manifolds $(M, \bold T)$ with the property that each edge has
distinct end points.

\medskip
\section{\bf Circle valued angle structures and maximization of volume}

Following Casson, Rivin \cite{Ri1} and Lackenby \cite{Lac}, we
introduced the following notion in \cite{Lu1},

\medskip
\noindent {\bf Definition.}  An \it $\bold S^1$-angle structure
\rm (or SAS for simplicity) on a closed triangulated pseudo
3-manifold $(M, \bold T)$ is a function $x: \Box ( \bold T) \to
S^1$ so that

(5) for each tetrahedron $\sigma$, $\prod_{ q \subset \sigma} x(q)
= -1$; and

(6) for each edge $e \in E$,  $\prod_{ q \in \Box} x(q)^{ i(e, q)}
= 1$.

\medskip

Let $SAS(\bold T)$ be the set of all $\bold S^1$-angle structures
on the triangulation $\bold T$.  If $x \in SAS(\bold T)$ and $v
\in TAS(\bold T)$, then $xe^{iv}$,  defined by $xe^{iv}(q) = x(q)
e^{iv(q)}$, is still in $SAS(\bold T)$. We use this to identify
the tangent space of $SAS(\bold T)$ with $TAS(\bold T)$.
 The \it Lobachevsky-Milnor volume \rm
(or simply the \it volume\rm) of an $\bold S^1$-angle structure
$x$ is defined to be:
$$ vol(x) = \sum_{ q \in \Box} \Lambda( arg(x(q)))$$
where $arg(w)$ is the argument of a complex number $w$ and
$\Lambda(t) = -\int_{0}^t \ln |2 \sin(s)| ds$. The volume formula
is derived from the volume of an ideal hyperbolic tetrahedron. See
Milnor \cite{Mil}.
 It is well known that $\Lambda(t): \bold R \to
\bold R$ is a continuous function with period $\pi$.  Thus,
$$ vol: SAS(\bold T) \to \bold R$$
is a continuous function.  Our goal is to relate the critical
points of $vol$ with the topology and geometry of the 3-manifold.

Using volume maximization to find geometric structures  based on
angle structures for manifolds with cusps was introduced by Casson
\cite{Ca} and Rivin \cite{Ri2}. In a recent work \cite{Gu},
Gu\'eritaud used the tool to prove the existence of hyperbolic
metrics on the once-punctured torus bundle over the circle with
Anosov holonomy.  Our approach follows the same path in a more
general setting.

\medskip
\subsection{Existence of SAS and critical points of volume}

In \cite{LT1}, we proved a general theorem on the existence of
real-valued prescribed-curvature angle structures on a
triangulated pseudo 3-manifold.  One can check that the proof in
\cite{LT1} implies the following proposition. Also see \cite{Lu1}
for a proof.

\medskip
\noindent {\bf Proposition 4.1 } (\cite{LT1}, \cite{Lu1}) \it For
a closed triangulated pseudo 3-manifold $(M, \bold T)$, the space
$SAS(\bold T)$ is non-empty and is a smooth closed manifold of
dimension $\chi(M) + |T|$. In particular, the volume $vol:
SAS(\bold T) \to \bold R$ has a maximum point. \rm

\medskip

Following \cite{LT1}, we give a short proof of it for real-valued
angle structures on ideally triangulated 3-manifolds with torus
boundary (i.e., closed triangulated pseudo 3-manifolds so that
each vertex link is a tours). The main idea of the proof for the
general case is the same.

Suppose otherwise that such a manifold $(M, \bold T)$ does not
support a real-valued angle structure.  Consider the linear map
$h:  R^{\Box} \to \bold R^T \times \bold R^E$ so that
$h(x)(\sigma)=\sum_{q \subset \sigma} x(q)$ and $h(x)(e) =\sum_{q}
i(e, q) x(q)$.  Let
 $\alpha \in \bold R^T \times \bold R^E$ be
$\alpha(\sigma)=\pi$ and $\alpha(e)=2\pi$.  Then the assumption
that $(M, \bold T)$ does not support a real-valued angle structure
means $ \alpha \notin h(\bold R^{\Box})$. Therefore, there exists
a vector $f \in \bold R^T \times \bold R^E$ so that $f$ is
perpendicular to the image $h(\bold R^{\Box})$ and $(f, \alpha)
\neq 0$.  This means that
$$\frac{1}{\pi} (f, \alpha) = \sum_{\sigma} f(\sigma) + 2 \sum_{e} f(e) \neq 0$$ and
$$ (h(x), f) = (x, h^*(f)) =0,$$
for all $x \in \bold R^{\Box}$, i.e.,
$$ h^*(f)=0,$$ where $h^*$ is the transpose of $h$.

Since $h^*(f)(q) =\sum_{\sigma, q \subset \sigma} f(\sigma) +
\sum_{e} i(e, q) f(e)$,  it follows that if $e, e'$ are two
opposite edges in $\sigma$, then
$$  f(e) + f(e') = -f(\sigma).$$
In particular, the sum of the values of $f$ at opposite edges in
$\sigma$ is independent of the choice of the edge pair. By lemma
2.1 (b), we see that there is a map $g$ defined on the pairs $(v,
\sigma)$ with $v < \sigma$ so that

$$ f(e) = g(v, \sigma) + g(v', \sigma)$$ where $v, v' < e$.
By the same argument as the one we used in the proof of  theorem
3.3, we see that $g(v, \sigma)$ is independent of the choices of
tetrahedra, i.e., $g: V \to \bold R$ so that

$$ f(e) = \sum_{ v < e} g(v)$$ and
$$ f(\sigma) = -\sum_{v < \sigma} g(v).$$

Now we express $$\sum_{\sigma} f(\sigma) + 2\sum_{e} f(e)$$
$$=-\sum_{\sigma, v < \sigma} g(v) +  2\sum_{e, v < v} g(v)$$
$$=-\sum_{v \in V}g(v) [\sum_{ \sigma > v} 1 -  2 \sum_{ e > v} 1]$$
$$=\sum_{v \in V} g(v)  [-|\{\text{ triangles in lk(v)}\}| +  2|\{ \text{vertices in lk(v)}\}|] =0$$
The last equality is due to the fact that the number of vertices
of a triangulation of the torus is equal to half of the number of
triangles in the triangulation. Also, in the summations
$\sum_{\sigma >v} 1$ and $\sum_{e > v} 1$, we count $\sigma$ and
$e$ with multiplicities, i.e., if $\sigma$ (or $e$) has $k$
vertices which are $v$, then $\sigma$ (or $e$) is counted $k$
times in the sum.  This ends the proof for manifolds with torus
boundary.

Proposition 4.1 guarantees that critical points for the volume
function always exist.  Here the concept of critical point of the
non-smooth function $vol$ has to be clarified.  It can be shown
(\cite{Lu1}) that for any point $p \in SAS(\bold T)$ and any
tangent vector $v$ of $SAS(\bold T)$ at $p$, the limit $\lim_{t
\to 0} \frac{ vol(p e^{itv}) - vol(p)}{t}$ always exists as an
element in $[-\infty, \infty]$.  A point $p \in SAS(\bold T)$ is
called a critical point of the volume if the above limit is  $0$
for all tangent vectors $v$ at $p$.

The main focus of our research is to extract topological and
geometric information from critical points of the volume function
on $SAS(\bold T)$.

Pursuing in this direction,  we have proved the following
\cite{Lu1}.

\medskip
\noindent {\bf Theorem 4.2.}  \it Let $(M, \bold T)$ be an
oriented triangulated closed pseudo 3-manifold. Suppose $x$ is a
critical point of the volume function $vol$ on $SAS(\bold T)$.

(a) If the critical point $x$ is a non-smooth point for the volume
function, so $x(q') =\pm 1$ for some $q' \in \Box$, then there
exists a solution $y$ to Haken's normal surface equation defined
on $\bold T$, which has exactly one or two non-zero quadrilateral
coordinates with $y(q') \neq 0$;

(b) If the critical point $x$ is a smooth point (i.e., $x(q) \neq
\pm 1$ for all $q$), then $x$ produces a solution to the
generalized Thurston equation. \rm

\medskip

Here are the key steps in the proof of theorem 4.2. Given $x \in
SAS(\bold T)$, we say a tetrahedron $\sigma \in T$ is \it flat \rm
with respect to $x$ if $x(q) = \pm 1$ for all $q \subset \sigma$
and \it partially flat \rm if $x(q) = \pm 1$ for one $q \subset
\sigma$. Let $U$ be the set of all partially flat but not flat
tetrahedra and $\bold W =\{ q | x(q) = \pm 1, q \subset \sigma,
\sigma \in U\}$.
 By analyzing the derivative of
 $\int_0^t \ln|2 \sin(s)| ds$,  we obtain the
following main identity. For $u \in TAS(\bold T)$,

\begin{equation}
\frac{d}{dt} vol( x e^{iut}) =-\sum_{q \in \bold W} u(q) \ln|t| -
\sum_{ q, x(q) =\pm 1} u(q) \ln|u(q)| -\sum_{q, x(q) \neq \pm 1}
u(q) \ln | \sin ( arg( x(q)))| + o(t)
\end{equation}

Now at a smooth point $x$, equation (4.1) becomes

$$ \frac{d}{dt} vol( x e^{iu})|_{t=0} = -\sum_{q}
u(q) \ln | \sin ( arg( x(q)))|.$$

By taking $u = u_e$ given by (3.10) and assuming $x$ is a smooth
critical point, we obtain a solution $z \in \bold C^{\Box}$ to the
generalized Thurston equation, where

$$  z(q) = \frac{ \sin( arg(x(q')))}{\sin( arg(x(q'')))} x(q),$$
and $ q' \to q \to q''$.  This argument was known to Casson
\cite{Ca} and Rivin. One may find a detailed argument in
\cite{Lu1} or \cite{Futer}.

If $x$ is a non-smooth critical point, then we deduce from (4.1)
two equations for all $u \in TAS(\bold T)$,

\begin{equation}
\sum_{q \in \bold W} u(q) =0, \end{equation}

and

\begin{equation}
\sum_{ q, x(q) =\pm 1} u(q) \ln | u(q)| =g(u),
\end{equation}
where $g(u)$ is a linear function in $u$.  Now we use the
following simple lemma.

\medskip\noindent {\bf Lemma 4.3.} \it Suppose $\bold V$ is a finite
dimensional vector space over $\bold R$ and $f_1, ..., f_n, g$ are
linear functions on $\bold V$ so that for all $x \in \bold V$,

 $$
\sum_{ i=1}^n f_i(x) \ln | f_i( x)| = g(x). $$

 Then for each index
$i$ there exists $ j \neq i$ and $\lambda_{ij} \in \bold R$ so
that
$$f_i(x) = \lambda_{ij} f_j(x).$$ \rm

\medskip
Using lemma 4.3 for (4.3) where the vector space $\bold V$ is
$TAS(\bold T)$ and the linear functions are $ u(q)$ with $x(q)
=\pm 1$ and $g$, we conclude that for each $q$ with $x(q) = \pm
1$, there exist $q_1$ and $\lambda \in \bold R$ so that $ u(q)=
\lambda u(q_1)$ for all $u$ $\in TAS(\bold T)$.
  This shows that for all $u
\in TAS(\bold T)$, the inner product $( u, q^* -\lambda q_1^*)=0$.
By theorem 2.2,  $q^* -\lambda q_1^*$ is in $Proj_{\Box}(\bold
S_{ns})$. Thus theorem 4.1 (a) follows.

\subsection{ Futer-Gu\'eritaud's Theorem}

In an unpublished work \cite{FG}, David Futer and Francois
Gu\'eritaud proved a very nice theorem concerning the non-smooth
maximum points of the volume function. The proof given below is
supplied by Futer and Gu\'eritaud. We are grateful to Futer and
Gu\'eritaud for allowing us to present their proof in this paper.

\medskip
\noindent {\bf Theorem 4.4} (Futer-Gu\'eritaud) \it Suppose $(M,
\bold T)$ is an oriented triangulated closed pseudo 3-manifold. If
$x$ is a non-smooth maximum point of the volume function on
$SAS(\bold T)$, then there exists a non-smooth maximum volume
point $y \in SAS(\bold T)$ so that all partially flat tetrahedra
in $y$ are flat. \rm

\medskip
\noindent {\bf Proof}(Futer-Gu\'eritaud).   Suppose $x$ is a
non-smooth maximum volume point in $SAS(\bold T)$.   Let $U$ be
the set of all partially flat but not flat tetrahedra in $x$ and
$\bold W =\{ q \subset \sigma | \sigma \in U$, $ x(q) = \pm 1$\}
as above. Note that, by assumption, for each tetrahedron $\sigma$,
there is at most one quadrilateral in $\bold W$ contained in
$\sigma$.

\medskip
\noindent {\bf Claim.}  Define the vector $ v =\sum_{ q' \in \bold
W} \sum_{ q \in \Box} w(q', q) q^*$, i.e., $v(q) = \sum_{ q' \in
\bold W} w(q', q)$. Then  $v \in TAS(\bold T)$.

\medskip

To see the claim, we must verify two conditions for $v$: (1) for
tetrahedron $\sigma \in \bold T$, $\sum_{ q \subset \sigma} v(q)
=0$, and (2) for each edge $e \in E$, $\sum_{q' \in \Box} i(e, q')
v(q') =0$.

The first condition follows from the fact that for any $q \in
\Box$, $\sum_{ q' \in \Box} w(q', q)=0$.  Indeed, for each
tetrahedron $\sigma$, $\sum_{ q \subset \sigma} v(q) = \sum_{q'
\in \bold W} \sum_{q \subset \sigma} w(q', q) =0$.

To see the second condition, we use (4.2).

By (3.9), for $e \in E$, the vector
$$ u_e= \sum_{e} W(e, q) q^*,$$ i.e., $u_e(q) = W(e, q)$, is in $TAS(\bold T)$.
 Taking this $u_e$ to be the vector $u$ in (4.2), we obtain

$$ \sum_{ q \in \bold W} W(e, q)
=0$$ i.e.,

$$\sum_{ q' \in \Box} i(e, q') (\sum_{q \in \bold W} w(q', q) )
=0.$$   The last equation says $\sum_{q'} i(e, q') v(q') =0$. This
verifies the claim.

\medskip

Now back to the proof of the theorem. For each point $p \in
SAS(\bold T)$, let $N(p)$ be the number of partially flat but not
flat tetrahedra in $p$. For the maximum point $x$, we may assume
$N(x)
>0$. We will produce a new maximum point $y$ so
that $N(y) < N(x)$ as follows. Let $v$ be the tangential angle
structure constructed in the claim above. Consider the smooth path
$$ r(t) = x e^{itv} \in SAS(\bold T).$$
Note, by definition, for $|t|$ small, $N(r(t)) =N(x)$. Take
$|t_0|$ be the smallest number so that
$$N(r(t)) =N(x) $$ for all $|t| < |t_0|$ and
$$N(r(t_0)) < N(x).$$
Furthermore, due to the basic property of the Lobachevsky function
that
$$ \Lambda(a) + \Lambda(b) + \Lambda(c) =0$$
for $a+b+c \in  \pi \bold Z$ and one of $a,b, c$ is in $\pi \bold
Z$, we have $$ vol(r(t)) = vol(x)$$ for $|t| \leq |t_0|$. Take
$y=r(t_0)$. Then we have produced a new maximum point $y$ with
smaller $N(y)$. Inductively, we produce a new maximum point so
that all partially flat tetrahedra are flat.  This ends the proof.

\medskip
\noindent

Combining theorem 4.2 with the theorem of Futer-Gu\'eritaud, we
obtain theorem 1.2,

\medskip
\noindent {\bf Theorem 1.2} \it Suppose $(M, \bold T)$ is a closed
triangulated oriented pseudo 3-manifold. Then there either exists
a solution to the generalized Thurston equation or there exists a
cluster of three 2-quad-type solutions to Haken's normal surface
equation.  \rm

\medskip

Indeed, by Futer-Gu\'eritaud's theorem, we can produce a
non-smooth maximum point $y$ so that there are three distinct
quadrilaterals $q_1, q_2, q_3$ in a tetrahedron with $y(q_i) = \pm
1$. Now we use theorem 4.2 to produce the corresponding
2-quad-type solutions $x_i$, one for each $q_i$ with $x_i(q_i)
\neq 0$. Note that we do not assume that $x_1, x_2, x_3$ are
pairwise distinct.

A stronger version of conjecture 1 is the following.

\medskip
\noindent{\bf Conjecture 2}  Suppose $(M, \bold T)$ is a minimally
triangulated closed irreducible 3-manifold so that all maximum
points of the volume function $vol: SAS(\bold T) \to \bold R$ are
smooth for $vol$. Then Thurston's equation on $\bold T$ has a
solution.

\subsection{\bf  Minimal triangulations with a cluster of three 2-quad-type solutions}

Our recent joint work with Stephan Tillmann shows the following.

\bigskip
\noindent {\bf Theorem 1.3} (\cite{LT2})  \it Suppose $(M, T)$ is
a minimally triangulated orientable closed 3-manifold so that
there are three 2-quad-type solutions $x_1, x_2, x_3$ of Haken's
equation with $x_i(q_i) \neq 0$ for three distinct quadrilaterals
$q_1, q_2, q_3$ inside a tetrahedron. Then,

(a) $M$ is reducible, or

(b) $M$ is toroidal, or

(c) $M$ is a Seifert fibered space, or

(d) $M$ contains the connected sum $\#_{i=1}^3 RP^2$ of three
copies of the projective plane. \rm

\medskip
By the work of W. Thurston and others, it is known, without using
the Ricci flow method,  that manifolds in class (d) but not in
cases (a), (b), (c) above are either Haken or hyperbolic. See for
instance \cite{smith}. Indeed, an irreducible, non-Haken,
atoroidal, non-Seifert-fibered 3-manifold containing  $\#_{i=1}^3
RP^2$ has a two fold cover which is a closed 3-manifold of
Heegaard genus at most 2. Such a manifold admits a $\bold Z_2$
action with 1-dimensional fixed point set. By Thurston's Orbifold
theorem \cite{BLP}, or \cite{CHK} one concludes that the manifold
is hyperbolic.

The proof of theorem 1.3 makes essential uses of Jaco-Rubinstein's
work on 0-efficient triangulations.
 We analyze carefully the cluster of three  2-quad-type solutions of Haken's
normal surface equation constructed from theorem 1.2.

Theorem 1.3 takes care of the topology of closed minimally
triangulated 3-manifolds which have non-smooth maximum volume
points.

We don't know if theorem 1.3 can be improved by using only one
2-quad-type solution instead of a cluster of three 2-quad-type
solutions.  Such an improvement will help in reproving the
Poincar\'e conjecture. For instance, one can weaken conjecture 1
by replacing the cluster of three 2-quad-type solutions by one
2-quad-type solution.  Another related conjecture is the
following,

\medskip
\noindent {\bf Conjecture 3} Suppose $(M, \bold T)$ is a minimally
triangulated closed orientable 3-manifold so that one edge of
$\bold T$ has the same end points and is null homotopic in $M$.
Then there exists a cluster of three 2-quad-type solutions on
$\bold T$.

\medskip
By theorem 1.3, one sees that conjecture 3 implies the Poincar\'e
conjecture without using the Ricci flow.

\medskip

\section{Some open problems}

Another potential approach to conjecture 1 is to use volume
optimization on a space closely related to $SAS(\bold T)$.  Let
$W(\bold T)$ be the space $\{  z \in \bold C^{\Box} | $ so that if
$ q \to q'$, then $z(q') = 1/(1-z(q))$ and for each edge $e$, the
right-hand-side of (3.7) is a positive real number\}.  The volume
function $vol: W(\bold T) \to \bold R$ is still defined. The
maximum points of the volume are related to the solutions of
Thurston's equation. In fact, a critical point of the volume
function in the set $W(\bold T) \cap (\bold C - \bold R)^{\Box}$
gives a solution to Thurston's equation.

It is conceivable that the following holds.

\medskip
\noindent {\bf Conjecture 4} Suppose $(M, \bold T)$ is a closed
orientable triangulated 3-manifold so that $W(\bold T) \neq
\emptyset$. Then $\sup\{ vol(z) | z \in W(\bold T)\} \leq v_3
||M||$ where $v_3$ is the volume of the regular ideal hyperbolic
tetrahedron.

\medskip
The first step to carry out this approach is to find conditions on
the triangulation $\bold T$ so that $W(\bold T)$ is non-empty.  To
this end, we consider solving Thurston's equation over the real
numbers, i.e., $z \in \bold R^{\Box}$.
 Here is a step toward
producing a real-valued solution to Thurston's equation.

\medskip
\noindent {\bf Definition 5.1}. Let $\bold Z_2$ be the field of
two elements $\{0,1\}$. A \it $\bold Z_2$-taut structure \rm on a
triangulated closed pseudo 3-manifold $(M, \bold T)$ is a map $f:
\Box \to \{0,1\}$ so that

(a) if $q_1, q_2, q_3$ are three quadrilaterals in each
tetrahedron $\sigma$, then exactly one of $f(q_1), f(q_2), f(q_3)$
is 1, and

(b) for each edge $e$ in $\bold T$, $\sum_{ q \in \Box} i(e,q)
f(q) =0$.

\medskip

The motivation for the definition comes from taut triangulations
and real-valued solutions to Thurston's equation. Indeed, if $z$
is a real-valued solution to Thurston's equation, then there is an
associated $\bold Z_2$-taut structure $f$ defined by:  $f(q) =0$
if $z(q) >0$ and $f(q)=1$ if $z(q) < 0$.  Another motivation comes
from taut triangulations. Suppose $\bold T$ is a taut
triangulation, i.e., there is a map $g: \Box \to \{0, \pi\}$ so
that for each tetrahedron $\sigma$, $\sum_{ q \subset \sigma} g(q)
= \pi$ and for each edge $e$, $\sum_{q} i(e,q) g(q) = 2\pi$. Then
one defines a $\bold Z_2$-taut structure by $f(q) = \frac{1}{\pi}
g(q)$.  A very interesting question is to find condition on $\bold
T$ so that $\bold Z_2$-taut structures exist. Is it possible that
the non-existence of $\bold Z_2$-taut structures implies the
existence of some special solutions to Haken's normal surface
equation?

Tillmann and I observed that the equations for $\bold Z_2$-taut
structures are non-linear but quadratic in $f(q)$. Indeed, a
vector $f \in \bold Z_2^{\Box}$ is a $\bold Z_2$-taut structure if
and only if condition (b) in definition 5.1 holds and for each
tetrahedron $\sigma$
\begin{equation}
 \sum_{ q \subset \sigma} f(q) =1, \end{equation} and
$$\sum_{q \neq q', q,q' \subset \sigma} f(q) f(q') =0.$$

The condition (b) in definition 5.1 and (5.1) should be considered
as the definition of a $\bold Z_2$-angle structure.

We end the paper with  several questions.

\noindent{\bf Question 1.} Given a triangulated pseudo 3-manifold
$(M, \bold T)$, when does there exist a $\bold Z_2$-taut
structure?  Can one relate the non-existence of $\bold Z_2$-taut
structure to some special solutions to Haken's equation?

\noindent {\bf Question 2.} When is a critical point of the volume
function of Morse type (i.e., when is the Hessian matrix
non-degenerated) and  when is the volume function  a Morse
function?

Let $v_3$ be the volume of the ideal regular hyperbolic
tetrahedron.

\noindent {\bf Question 3.} Is the Gromov norm of a closed
3-manifold multiplied by $v_3$ among the critical values of the
volume function?

\noindent {\bf Question 4.} Is it possible to produce a Floer type
homology theory associated to the volume function on $SAS(\bold
T)$ which will be a topological invariant of the 3-manifold?

Department of Mathematics

Rutgers University

New Brunswick, NJ 08854, USA

email: fluo\@math.rutgers.edu
\end{document}